
\tolerance=10000
\hsize=13truecm
\hoffset=1.4truecm

\expandafter\ifx\csname amssym.def\endcsname\relax \else\endinput\fi
\expandafter\edef\csname amssym.def\endcsname{%
       \catcode`\noexpand\@=\the\catcode`\@\space}
\catcode`\@=11
\font\tenmsb=msbm10
\font\sevenmsb=msbm7
\font\fivemsb=msbm5
\newfam\msbfam
\textfont\msbfam=\tenmsb
\scriptfont\msbfam=\sevenmsb
\scriptscriptfont\msbfam=\fivemsb
\def\Bbb#1{{\fam\msbfam\relax#1}}


\def\Z{{\Bbb Z}}



\def\H{{\cal H}}

\def\AA{{\cal A}}
\def\AI{\AA ({I_0})}

\def\AP{{\AI^+}}
\def\BBB{{\cal B}}

\def\ss{\sigma}
\def\th{{\theta}}



\def\sm{\setminus}
\def\inn{\subseteq}
\def\nni{\supseteq}

\def\bull{\hfill\vrule height 1.6ex width .7ex depth -.1ex }
\def\vid{\emptyset}
\def\bs{\backslash}

\def\Id{{\rm Id}}
\def\Hom{{\rm Hom}}

\def\Homgr{{\rm Homgr}}

\def\Cb{{\bf C}^{\bf b}}
\def\Kb{{\bf K}^{\bf b}}
\def\Db{{\bf D}^{\bf b}}

\def\mod{{\bf mod}}

\def\te{\otimes}
\def\ddd{\partial}

\def\sg#1{(-1)^{#1}}
\def\mr#1{\smash{\mathop{\relbar\joinrel\longrightarrow}\limits^{#1}}}

\centerline{\bf On Okuyama's theorems about Alvis-Curtis duality.
} \bigskip

\centerline{Marc Cabanes
} 
\bigskip

\noindent{\bf Abstract : } {\it We report on theorems by T. Okuyama about complexes
generalizing the Coxeter complex and the action of parabolic subgroups on them, both for
finite
 BN-pairs and finite dimensional Hecke algebras. Several simplifications,
using mainly the surjections of [CaRi], allow a slightly more compact treatment than the
one in [O].}

\footnote{}{2000 \it Mathematics Subject Classification. 20C08, 20C20, 20C33, 20J05.}
\bigskip

The purpose of this paper is to report on the unpublished manuscript [O] by T. Okuyama
where are proved  some conjectures generalizing to homotopy categories the theorems
of [CaRi] and [LS] holding in derived categories. We refer to the latter references and
[CaEn]~\S 4 for a broader introduction to the subject. 

The main theme is the one of complexes related with the Coxeter complex and the action
of parabolic subgroups on them, either for finite groups with BN-pairs or for finite
dimensional Hecke algebras. Okuyama's contractions prove a quite efficient tool in a
number of situations (see the proof of Solomon-Tits theorem in
\S~6).

We often stray away from Okuyama's proofs when it allows simplifications. 

We also
emphasize some statements that may be of independent interest (see \S~1 below), and
actually re-prove [Ri1]~\S 8, [LS], and [CaRi].

Most proofs are selfcontained apart from basic facts about split
BN-pairs and Hecke algebras.\medskip

\noindent{\bf Notations.} Let $A$ be a ring, one denotes by $\mod(A)$ the associated
category of finitely generated left $A$-modules. If $B$ is another ring, $A,B$-bimodules
are just objects of $\mod(A\te B^\circ )$ where  $B^\circ $ denotes the opposite algebra
and $\te$ denotes (commutative) tensor product over $\Z$. Whenever $M$ is an
$A,B$-bimodule and $N$ is a $B$-module, one denotes by $M\te_BN$ the usual tensor
product over $B$, considered as an $A$-module. When $B$ is a group algebra $RG$ with $R$
a commutative ring and $G$ a finite group, one may use the abbreviation $M\te_GN$ for
$M\te_{RG}N$.

We denote by
$\Cb (A)$,
$\Kb(A )$, and $\Db(A)$ the categories of
bounded (cochain) complexes of
$A$-modules, its homotopy and derived categories, respectively. We allow to view
a given complex of $A$-modules $(\dots C^i\mr{\ddd ^i}  C^{i+1}\dots )$ as a
$\Z$-graded
$A$-module
$C:=\oplus_i C^i$ endowed with a homogeneous endomorphism $\ddd_C=(\ddd^i)_i$
of degree 1 satisfying
$\ddd_C\circ \ddd_C =0$.

Let
$S$ be a finite set, we denote by
$2^S$ the set of subsets of
$S$. 
A {\it coefficient system} on $2^S$ is a family of
$A$-modules $(M^I)_{I\inn S}$ and $A$-homomorphisms $\varphi_{JI}^M\colon M^I\to M^J$
defined when $I\inn J$ (``restriction maps") and satisfying
$\varphi_{KJ}^M\circ\varphi_{JI}^M=\varphi_{KI}^M$ whenever $I\inn J\inn K\inn S$.

Choose a total ordering on $S$. When $I\inn S$ and $s\in S$,
denote by
$n(I,s)$ the number of elements of $I$ that are $< s$. A coefficient system
$((M^I)_I, (\varphi_{JI}^M)_{I,J})$ gives rise to an object $(\dots\to M^0\to
M^1\to\dots )$ of
$\Cb(A )$, where $M^i=\oplus_{I;|I|=i}M^I$ and
$\ddd^i\colon M^i\to M^{i+1}$ is defined on
$M^I$ by $\ddd^i{}_{|M^I}=\sum_{s\in S\sm I}(-1)^{n(I,s)}\varphi_{I\cup \{ s\},I}^M$. For a
more canonical definition, showing independence with regard to the choice of an ordering,
 see for instance
[CaEn]~Exercise~4.1.

\bigskip\centerline {\bf 1. Reduced elements in Coxeter groups and Okuyama's
contractions.}

\medskip
Assume $(W,S)$ is a finite Coxeter group. One denotes by $l$ the length map
with respect to $S$. 
If $I,J\inn S$, denote by $W_I$ the subgroup of $W$ generated by $I$ and by $D_{IJ}$
the set of $w\in W$ such that $l(uwv)=l(u)+l(w)+l(v)$ for any $u\in W_I$, $v\in W_J$.

If $V,V'$ are subsets of $W$, one denotes $VV'=\{ vv'\mid v\in V,\ v'\in V'\}$ and
$V^{-1}=\{ v^{-1}\mid v\in V\}$.

Let
$\AA =\cup_{I\inn S}W_I\backslash W=\{ W_Iw\mid I\inn S,\ w\in W\}$ be the set of
right cosets mod. parabolic subgroups. 

If $a\in\AA$, denote $S(a):=S\cap aa^{-1}$ (that is $S(a)=I$ if $a=W_Iw$ for
some $w\in W$) and
$l(a)={\rm min}\{ l(x)\mid x\in a\}$, which is the length of the only element of
$a\cap D_{S(a),\vid}$.

If $J\inn S$ and $a\in \AA$, denote $a\cup J:=W_{J\cup S(a)}a$. If $s\in S$,
we write $a\cup s=a\cup  \{ s\} $.

The {\sl Coxeter complex} $A(W,S)$ is the complex associated to the coefficient
system on $2^S$ defined by $I\mapsto \Z (W_I\bs W)=\oplus_{a\in \AA ,\ S(a)=I}\Z
a$ with restriction maps $\varphi_{JI}$ defined by $\varphi_{JI}
(a)=a\cup J$ for $I\inn J$. Then $A(W,S)=\Z\AA$ graded by $a\mapsto |S(a)|$ and
the differential is defined by $\ddd (a)=\sum_{s\in S\sm S(a)}
(-1)^{n(S(a),s)}a\cup s$ for $a\in\AA$ (recall that $S$ is ordered and $n(I,s )$
denotes the number of elements of $I$ that are $< s$).

Fix $I_0\inn S$, $I_0\not= S$ (the case $I_0=S$ is trivial for what follows but
could create ambiguities).

\medskip\noindent {\bf Notation~1. } Let $\AI :=\{ a\in \AA\mid a\cap D_{\vid
,I_0}\not=\vid\}$. Note that the condition defining $\AI $ is equivalent to the element
of minimal length in $a$ being in $D_{\vid
,I_0}$, hence in $D_{S(a)
,I_0}$. Let
$\AP :=\{ a\in\AI\mid a\not\inn W_{I_0}\}$, that is the classes $W_Iw$ for $I\inn S$, $w\in
D_{I,I_0}$, and $w\not= 1$ if moreover $I\inn I_0$.
\hfill\break If
$a\in\AI$, denote $I_0(a)= I_0\cap a^{-1}a$. 

\medskip

Note the following property

\noindent{\bf (P2)} If $a\in \AI$, then $I_0(a)=I_0\cap I^w$ for $I=S(a)$ and the unique
$w\in a\cap D_{I,I_0}$ (see
[CaEn]~2.6).

\medskip\noindent {\bf Proposition~3.} {\sl  $\Z\AP$ is a subcomplex of $A(W,S)$
in
$\Cb (\Z )$. There exists a linear map
$\ss\colon  \Z\AP\to \Z\AP$ homogeneous of degree $-1$ such that 

(i) $\ss \ddd+\ddd\ss=\Id$.

(ii)  for any $a\in \AP$, $\ss (a)\in\oplus \Z b$ where the sum is over
$b\in\AP$ such that $I_0(b)\nni I_0(a)$.
}

\medskip

The Proposition will be used in \S~2-4 through the following.

\medskip\noindent {\bf Theorem~4.} {\sl Keep $(W,S)$ a finite Coxeter
group, and
$I_0\inn S$.  Let $A$ be a ring and
$((M^I)_I, (\varphi_{JI}^M)_{I,J})$ be a coefficient system of $A$-modules
on $2^{I_0}$.

Assume that for all $b\in \AP$, we are given a submodule $Z_b\inn
M^{I_0(b)}$ such that, for all $b, b'\in \AP$ with $b\inn b'$ (and therefore $I_0(b)\inn
I_0(b')$), one has
$\varphi^M_{I_0(b'),I_0(b)}(Z_b)\inn Z_{b'}$ .

Let a coefficient system $Z$ on $2^S$ be defined by
$I\mapsto Z^I:=\Pi_{b\in\AP ; S(b)=I}Z_b$, and
$\varphi^Z_{I',I}(x)=\varphi^M_{I_0(b\cup I'),I_0(b)}(x)\in Z_{b\cup I'}\inn
Z^{I'}$ if
$I\inn I'\inn S$ and
$x\in Z_{b}$ with $S(b)=I$. 

Then the complex associated
to $Z$ is contractible in $\Cb (A)$. }

\medskip\noindent {\it Proof.} 
Let $\ss\colon
\Z\AP\to \Z\AP$ be a map as in Proposition~3. This implies the existence of integers
$m_{b,b'}$ for
$b,b'\in\AP$ with ${I_0}(b)\inn {I_0}(b')$ and $|S(b')|=|S(b)|-1$ such that $\ss
(b)=\sum_{b'}m_{b,b'}b'$. Set $m_{b,b'}=0$
for pairs $(b,b')$ not satisfying the conditions
above.

Since the differential on
$\Z\AP$ is defined by $\ddd (b)=\sum_{s\in S\sm
S(b)}(-1)^{n(S(b),s)}(b\cup s)$, the relation in
Proposition~3.(i) reads, for any $b\in\AP$,

$$ b=\sum_{\scriptstyle s\in S\sm S(b)\atop \scriptstyle b'\in\AP
}(-1)^{n(S(b),s)}m_{b\cup s,b'}
b'+ \sum_{\scriptstyle b'\in\AP\atop \scriptstyle s\in S\sm
S(b')}(-1)^{n(S(b'),s)}m_{b,b'}.(b'\cup s)\leqno({\bf
E5}_b)$$ in $\Z\AP$.

When $b\in\AP$ and $z\in Z_b$, we write $z_b$ to mean the element of
$Z^{S(b)}$ in the factor
$Z_b $.
Define $\bar\ss$ on $Z_b\inn
\prod_{b\in\AP}Z_b$  by
$\bar\ss (z)=\sum_{b'\in\AP}m_{b,b'}\varphi^M_{{I_0}(b'),{I_0}(b)}(z)_{b'}$.
This is well defined since $m_{b,b'}\not= 0$ implies
$I_0(b)\inn I_0(b')$.

The differential on $Z$ is defined by
$$\partial^{Z} (z)=\sum_{s\in S\sm S(b)}(-1)^{n(S(b),s)}
\varphi^M_{{I_0}(b\cup s),{I_0}(b)}(z)_{b\cup s}$$ for $z\in Z_b$ (remember
that
$Z_b$ is at degree
$|S(b)|$).  Then, keeping $z\in Z_b$,
$(\bar\ss\partial^Z +\partial^Z\bar\ss )(z) =\bar\ss
(\sum_{s\in S\sm S(b)}(-1)^{n(S(b),s)} \varphi^M_{{I_0}(b\cup
s),{I_0}(b)}(z)_{b\cup s}) + \partial^Z
(\sum_{b'\in\AP}m_{b,b'} \varphi^M_{{I_0}(b'),{I_0}(b)}(z)_{b'})$ which equals 

\leftline{$
\sum_{s\in S\sm S(b), b'\in\AP }(-1)^{n(S(b),s)}m_{b\cup s,b'}
\varphi^M_{{I_0}(b'),{I_0}(b)}(z)_{b'} +$}
\rightline{$ \sum_{b'\in\AP, s\in S\sm
S(b')} (-1)^{n(S(b'),s)} m_{b,b'} \varphi^M_{{I_0}(b'\cup
s),{I_0}(b)}(z)_{b'\cup s} $} 

\noindent thanks to the composition of
restriction maps. By {\bf
E5}$_b$ above, this equals $\varphi^M_{{I_0}(b),{I_0}(b)} (z)_b = (z)_b$.
Then
$\bar\ss\ddd +\ddd\bar\ss = \Id$, that is the contractibility of
$Z$ in $\Cb (A)$.
\bull

\medskip\noindent {\it Proof of Proposition~3.} 
That $\Z\AI$ and $\Z\AP$ are sub-coefficient systems of $\Z\AA$ is clear
by the definition of restriction maps and the fact that both $\AI$ and
$\AP$ are closed for supsets.

In order to check the particular contractibility that is announced about $\Z\AP$,
it seems handy to apply the following symmetry of $D_{\vid ,I_0}$ and $\AI$.

Let $\th \colon W\to W$ defined by $\th (w)=w_Sww_{I_0}$, where
$w_I\in W_I$ denotes the element of maximal length for $I\inn S$. This
$\th$ clearly preserves
$\AA$ with $S(\th (a))=S(a)^{w_S}$ and $\varphi_{JI}\circ\th =
\th\circ\varphi_{J^{w_S},I^{w_S}}$ for any $a\in\AA$ and $I\inn J\inn S$. 
Moreover $\th$ preserves
$D_{\vid, I_0}
$ (exchanging
$1$ and
$w_Sw_{I_0}$) and therefore
$\AI$, with $\th (D_{I,I_0})=D_{I^{w_S},I_0}$. One has also $I_0(\th(a))=I_0(a)^{w_{I_0}}$ for
all
$a\in
\AI$.

Denote $\BBB =\th (\AP )=\{ W_Iw\mid I\inn S, \ w\in D_{I,I_0},\
w\not=w_Sw_{I_0}\}$.

Note that if $w\in D_{\vid ,I_0}\sm\{ w_Sw_{I_0}\}$, there is some $s_w\in S\sm
(S\cap   wI_0w^{-1})$ such that
$l(s_ww)>l(w)$ (take the $w_S$-conjugate of the first term in any reduced
decomposition of $\th (w)\not= 1$). 

If $b\in\BBB$, $b=W_Iw$ with $I\inn S$ and
$w\in D_{II_0}$, one denotes
$s_b:=s_w$ and $b\sm s_b = W_{I\sm \{ s_b\}}w$.

\medskip\noindent {\bf Lemma~6.} {\sl 
 Let
$\tau \colon  \Z\BBB\to \Z\BBB$ be the linear map homogeneous of degree $-1$
defined as follows on
$b\in\BBB $
\smallskip\centerline{$\tau ( b) = (-1)^{n(S(b),s_b)} (b\sm s_b)$ if $s_b\in S(b)$,
\ \ \ \ $\tau ( b) = 0$ otherwise.}
\smallskip\noindent One has $\tau^2=0$ and \hfill\break\noindent (i) If $\tau
(b)\not= 0$ then $\tau (b)=\pm b'$ where $b'\in \BBB $, $l(b')=l(b)$ and ${I_0}(b')=
{I_0}(b)$.\hfill\break\noindent
(ii) $(\tau \ddd +\ddd \tau  )(b)\in b+
\oplus_{b'} \Z b'$ the sum being over $b'\in \BBB $ such that
$l(b')<l(b)$ and
${I_0}(b')\nni {I_0}(b)$.}
\medskip

Let us show how this implies Proposition~3.

By (i) of Lemma~6, $\tau (\BBB )\inn \Z\BBB$.  By (ii), $\tau \ddd +\ddd \tau = \Id +\rho$
where $\rho$ is nilpotent and of degree 0. Now define
\hfill\break\noindent \centerline{$\tau ':=\tau -\tau \rho +\tau \rho ^2-\tau \rho
^3+\dots =
\tau (\tau \ddd +\ddd \tau )^{-1}$}\break\noindent
\noindent and let us check that it satisfies the conditions (i) and (ii) of Proposition~3.

First  $\tau '
\ddd +\ddd
\tau ' =\Id$ since $\ddd$ clearly commutes
with
$\tau \ddd +\ddd \tau $, hence with its inverse. 

As for (ii), that is $\tau ' (b)\in \oplus_{b'} \Z b'$ where the sum is over
$b'\in \BBB $ such that ${I_0}(b')\nni {I_0}(b)$, this is a
consequence of
$\tau '=\tau -\tau \rho +\tau \rho ^2-\dots $ with $\tau $ and
$\rho $ having the corresponding property by the Lemma.

This implies now Proposition~3 by defining $\ss =\th\circ \tau '\circ \th$ thanks to
the elementary properties of $\th$ with regard to the restriction maps and
$a\mapsto I_0(a)$. Note however that $\th\circ \ddd\circ \th$ is not exactly
$\ddd$ but the same twisted by $w_S$-conjugacy due to the property of $\th$
with regard to restriction maps $\varphi_{JI}$. A
correction consists in
adapting the ordering on
$S$~: choose first the $w_S$-conjugate of the ordering implicit in Proposition~3.
\bull

\medskip\noindent {\it Proof of Lemma~6.} 
(i) Write $b=W_Iw$ with $w\in D_{I{I_0}}$, $w\not= w_Sw_{I_0}$. One must assume
$s_w\in I$. Then
$\tau (b)=\pm  W_{I\sm \{ s_w\} }w$ with $w\in D_{I\sm  \{ s_w\},{I_0}}$, $w\not=
w_Sw_{I_0}$ and
$I_0(b\sm s_b)=(I\sm  \{ s_w\})^w\cap {I_0}= I^w\cap {I_0}=I_0(b)$ by (P) and
the definition of
$s_w$. It is also clear that $s_{b\sm s_b}=s_w=s_b$, so $\tau^2(b)=0$.

(ii) Note that if $b,b'\in \AA$ and $b\inn b'$ (inclusion of cosets), then 
${I_0}(b)\inn {I_0}(b')$ and $l(b')\leq l(b)$. 

From (i) and the definition of $\ddd $ it is then clear that we don't have
to worry about sets ${I_0}(b')$'s. So we concentrate on lengths $l(b')$. Thanks to the
above and (i) just proved, on evaluating
$ (\ddd \tau +\tau \ddd )(b)$ we must check that only one term, the one producing
$b$, has length remaining equal to $l(b)$ upon applying $\tau $ and $\ddd $.

  First case~:
$ s_b\in S(b)$. Then $\tau (b)=\sg{n(b, s_b)} (b\sm  s_b)$ with $S(b\sm s_b)=S(b)\sm
\{ s_b\}$ and $((b\sm s_b)\cup s_b)=b$. Then
$$(\ddd \tau +\tau \ddd )(b)=b+\sum_{s\in S\sm S(b)}\sg{n(b,
s_b)+n(b\sm s_b,s)}((b\sm s_b)\cup s)+\sum_{s\in S\sm S(b)}\sg{n(b,s)}\tau (b\cup
s)$$
In the above sums, one must spot the terms with length $l(b)$ ($=l(w)$).  In the first
$\sum$, $b\sm s_b=W_{I\sm \{ s_b\} }w$ and one may have $l(W_{I\sm
\{ s_b\}\cup \{ s\}}w)=l(w)$ only if $l(sw)>l(w)$, thus producing $((b\sm
s_b)\cup s) =W_{I\sm \{ s_b\}\cup \{ s\} }w$ in the above sum. For the second
$\sum$, by (i), length $l(b)$ is kept only if $l(b\cup s)=l(b)$. This means again
$l(sw)>l(w)$, thus giving a term
$\tau (b\cup s)=\sg{n(b\cup s,s_b)}W_{I\sm \{ s_b\}\cup \{ s\} }w$. 

In all, the two $\sum$ contribute terms of length $l(b)$ by the sum over $s\in
S\sm S(b)$ with
$l(sw)>l(w)$ of the terms $(\sg{n(b,
s_b)+n(b\sm s_b,s)}+\sg{n(b,s)+n(b\cup s,s_b)})(b\sm s_b\cup s)$.
This is $0$ as can be seen by calculating the component on $b\cup s$ of
$\partial\circ\partial (b\sm s_b)$ , one finds
$\sg{n(b,
s_b)+n(b,s)}+\sg{n(b\sm s_b,s)+n(b\cup s,s_b)}=0$.

Second case~: $b=W_Iw$ with $w\in D_{I{I_0}}$ and
$s_b\not\in S(b)=I$. Then $\tau (b)=0$ and $b=\pm \tau (b\cup s_b)$ where $b\cup
s_b=W_{I\cup \{ s_w\} }w\in\BBB$ with $w\in D_{I\cup \{ s_w\},{I_0}}$ and $l(b\cup
s_b)=l(b)=l(w)$. To get our claim, it suffices to show that $ \ddd \tau +\tau \ddd
-\Id$ sends
$\tau (b\cup s_b)$ into $\oplus_{b';l(b')<l(b)}\Z b'$. We have $(\ddd \tau +\tau \ddd
-\Id)\tau (b\cup s_b)=(\tau \ddd \tau -\tau )(b\cup s_b)=\tau (\ddd \tau +\tau
\ddd -\Id)(b\cup s_b)$ and our claim follows from the first case and
Lemma~6.(i). \bull

\medskip\noindent {\bf Remark~7.} When $a\in\AI$, denote by $v_0(a)$ the (unique)
element of
$a\cap D_{S(a),I_0}$, so that $a=W_{S(a)}v_0(a)$ with $l(a)=l(v_0(a))$.

Denote by $\leq_{\bf r}$ the right divisibility in $W$, that is $w'\leq_{\bf r} w$ if and
only if
$w=w''w'$ with lengths adding.
Any inclusion $b\inn b'$ in $\AI$ clearly implies $v_0(b')\leq_{\bf r} v_0(b) $. 

A quick inspection shows that in the above proof some relation $l(b')<l(b)$ may occur only
when in addition 
$v_0(b') <_{\bf r} v_0(b)$. So the map $\tau\colon  \Z\BBB\to \Z\BBB$ of Lemma~6
satisfies
$(\tau \ddd +\ddd \tau )(b)\in b+
\oplus_{b'} \Z b'$ where the sum is over $b'\in \BBB $ such that
$v_0(b')<_{\bf r} v_0(b)$ (instead of just $l(b')<l(b)$).

So Proposition~3 holds with a map $\ss$ satisfying $\ss (a)\in\oplus_b\Z b$ where the
sum is over $b\in \AP$ such that $w_Sv_0(b) w_{I_0} \leq_{\bf r} w_Sv_0(a) w_{I_0}$ and
$I_0(b)\nni I_0(a)$.
\medskip

\bigskip

\centerline {\bf 2. A theorem of Curtis type in the homotopy category.}

\medskip

Let $G$ be a finite group endowed with a split BN-pair of
characteristic the prime $p$ (see [CaEn]~2.20). We have subgroups $N$, $B$, $T\inn
B\cap N$. The quotient $W:=N/T$ is a Coxeter group for the subset $S\inn
W$. When $I\inn S$, the associated parabolic subgroup $P_I=BW_IB$ is
a semi-direct product $U_I.L_I$ for $U_I$ the largest normal $p$-subgroup of $P_I$
and $L_I$ a group with a split BN-pair associated to the subgroups
$N\cap L_I$, $B\cap L_I$, $T$ and the Coxeter group $W_I$.

Let $R =\Z [p^{-1}]$ or any commutative ring where $p$ is invertible. 

\medskip\noindent {\bf Notation~8. } If $I\inn S$, denote $e_I=|U_I|^{-1}\sum_{u\in
U_I}u$, an idempotent in the group algebra $RG$. Define the coefficient system
$X(G)$ of $RG$-bimodules on
$2^S$ by 
$X(G)^I=RGe_I\te_{P_I}e_IRG$ and restriction maps
$\varphi_{JI}\colon X(G)^I\to X(G)^J$ defined by $x\te_{P_I} y\mapsto x\te_{P_J} y$
whenever
$x\in RGe_I$ and $y\in e_IRG$.

\medskip\noindent {\bf Theorem~9. } (Okuyama, [O]~3.1) {\sl  Let ${I_0}\inn S$. Then
$X(G)e_{I_0}\cong RGe_{I_0}\te_{L_{I_0}}X(L_{I_0})$ in $\Kb(RG\te (RL_{I_0}{})^\circ  )$
and
$e_{I_0}X(G)\cong X(L_{I_0})\te_{L_{I_0}}e_{I_0}RG$ in $\Kb(RL_{I_0}\te RG^\circ )$.}

\medskip

The proof consists in giving a description of the kernel of the  surjection $X(G)e_{I_0}\to
RGe_{I_0}\te_{L_{I_0}}X(L_{I_0})$ introduced in [CaRi]~3.5, allowing to apply
Theorem~4. One will use repeatedly the following (see [HL]~3.1,  [CaEn]~Ex~5 p 53).

\medskip\noindent {\bf Proposition~10. } {\sl $\bullet$ $e_Ie_J=e_Je_I=e_I$
when
$I\inn J\inn S$.

$\bullet$ If $I,J\inn S$, $w\in D_{IJ}$ and $n\in N$ with $nT=w$, then
$e_Ine_J=e_{I\cap ^{w}J}ne_{ J}=e_{I}ne_{I ^{w}\cap J}=e_{I\cap ^{w}J}ne_{I
^{w}\cap J}$ }.

\medskip\noindent {\bf Definition. } {If $I\inn S$, $w\in D_{I , {I_0}}$, let
\smallskip\centerline{$X_{I,w}:=
RGe_{I}\te_{P_{I}}e_{I}RP_IwRP_{I_0}e_{I_0}$}\smallskip
\noindent  and
\smallskip\centerline{$Y_{I,w}:= RGe_{I^w\cap {I_0}}\te_{P_{I^w\cap
{I_0}}}e_{I^w\cap {I_0}}RP_{I_0}=RGe_{I_0}\te_{L_{I_0}}X(L_{I_0})^{I^w\cap
{I_0}}$}\smallskip
\noindent   both
$RG\te RP_{I_0}^\circ$-modules.}
\medskip

In the following propositions, keep $I\inn S$ and $w\in D_{I{I_0}}$.

\medskip\noindent {\bf Proposition~11. } {\sl $X_{I,w}\cong Y_{I,w}$ as
$RG\te RP_{I_0}^\circ$-module by a map sending $x\te_{P_I} \dot w y $ to
$xe_{I\cap{}^w {I_0}}\dot we_{I^w\cap {I_0}}\te_{P_{I^w\cap {I_0}}}e_{I^w\cap
{I_0}} y$ for any $x\in RGe_I$, $y\in RP_{I_0}e_{I_0}$, $\dot w\in N$ such that $\dot
wT=w$ .}

\medskip\noindent {\bf Proposition~12. } {\sl   If $ J\inn S$, $w'\in D_{J ,{I_0}}$
and $I^w\cap {I_0}\inn J^{w'}\cap {I_0}$, there is a $RG\te
RP_{I_0}^\circ$-morphism
$Y_{I,w}\to Y_{J,w'}$ sending $x\te_{P_{I^w\cap {I_0}}}y$ to $x\te_{P_{J^{w'}\cap
{I_0}}}y$ for $x\in RGe_{I^w\cap {I_0}}$, $y\in e_{J^{w'}\cap {I_0}}RP_{I_0} $.

If $I\inn J\inn S$ and $w'\in D_{J{I_0}}\cap W_Jw$ (which ensures $I^w\cap {I_0}\inn
J^{w'}\cap {I_0}$)  the above morphism corresponds to the restriction map
$X(G)^Ie_{I_0}\to X(G)^Je_{I_0}$ of
$X(G)e_{I_0}$ through the isomorphism of  Proposition~11 and the decomposition
$X(G)^Ie_{I_0}=RGe_I\te_{P_I}e_IRGe_{I_0}=
\oplus_{P_IwP_{I_0}}RGe_{I}\te_{P_{I}}e_{I}RP_IwRP_{I_0}e_{I_0}$. }

\medskip

The above will allow to describe the kernel of the surjection $X(G)e_{I_0}\to
RGe_{I_0}\te_{L_{I_0}}X(L_{I_0})$ defined in [CaRi]~\S~6 (see also [CaEn]~4.10). 

The
following will be useful to deduce an isomorphism in the homotopy category.

\medskip\noindent {\bf Proposition~13.  } Let ${\cal A}$ be an abelian
category.  Let $0\to Z\to
Y\to Y'\to 0$ be an exact sequence in $\Cb ({\cal A} )$ which is split in each degree and
such that $Z$ is contractible. Then $Y\cong Y'$ in $\Kb ({\cal A})$.

\medskip\noindent {\it Proof of Proposition~13.} More generally, it is well known that in
a short exact sequence which is split in each degree, the third term is always homotopy
equivalent to the mapping cone of the monomorphism.  \bull

\medskip\noindent {\it Proof of Theorem~9. } Let us consider $Y$ the
complex associated to the coefficient system on $2^S$ defined by
$I\mapsto \oplus_{w\in D_{I{I_0}}}Y_{I,w}$ and restriction maps
$Y_{I,w}\to Y_{J,w'}$ as in Proposition~12 when $I\inn J$ and
$W_Jw=W_Jw'$, 0 otherwise.

Then  $X(G)e_{I_0}$ is isomorphic to $Y$ thanks to the
isomorphisms of Proposition~11. 

In order to apply Proposition~13, we show that there is surjection of
complexes
$Y\to RGe_{I_0}\te_{L_{I_0}}X(L_{I_0})$ in
$\Cb(RG\te ( RL_{I_0} )^\circ )$, split in each degree, and with
contractible kernel.

Note that $Y_{I,w}$ only depends on the class $W_Iw\in\AI$ (see
Notation~1). We then use the notation
$Y_b=RGe_{I_0}\te_{L_{I_0}}X(L_{I_0})^{{I_0}(b)}$ for
$b\in\AI$ (see (P2)). Note that $Y$ is a coefficient system on
$2^S$ defined by
$I\mapsto
\oplus_b Y_b$ where the sum is over $b\in\AI$ such that $S(b)=I$, and the
restriction map is defined by sending $Y_b$ to $Y_{b'}$ only when $b'\nni
b$ and is then
$RGe_{I_0}\te_{L_{I_0}}\varphi_{I_0(b'),I_0(b)}$ where $\varphi$ is the restriction
map of $X(L_{I_0})$.

Let us define $Y'=\oplus_bY_b$ where the sum is over $b\in\AI$ such that $b\inn
W_{I_0}$ and same restriction maps as in $Y$. Then $Y'\cong
RGe_{I_0}\te_{L_{I_0}}X(L_{I_0})$ since each $b\inn W_{I_0}$ is equal to $W_{S(b)}$
with $I_0(b)=S(b)\inn {I_0}$.

We have a surjective map of coefficient systems of $RG\te
RL_{I_0}^\circ$-modules on
$2^S$~:
$Y\to Y'\cong RGe_{I_0}\te_{L_{I_0}}X(L_{I_0})$ sending  $Y_b$ to $0$ whenever
$b\in\AI$ and
$b\not\inn W_{I_0}$. The kernel is $\oplus_bY_b$ where the sum is over $b\in\AP$.

One may now apply Theorem~4 with $A=RG\te
RL_{I_0}^\circ$, $M=RGe_{I_0}\te_{L_{I_0}}X(L_{I_0})$ and $Z_b=Y_b$ for
$b\in \AP$. This tells us that $\oplus_{b\in\AP}Y_b$ is a
contractible complex.

We then have the first isomorphism of the theorem. In
order to deduce the second one, one uses the (covariant) functor
$M\to M^\iota$ from $\mod(RG)$ to
$\mod(RG^\circ )$ which consists in keeping the same $R$-module structure and
composing action of group elements $g\in G$ with inversion (see [CaEn]~4.16 and
proof). It is exact, commutes with tensor products, extends to complexes, and
induces involutory equivalences of
$\Kb$ categories. It is easy to see that
$X(G)^\iota
\cong X(G)$,
$(RGe_{I_0})^\iota\cong e_{I_0}RG$ and $X(L_{I_0})^\iota\cong X(L_{I_0})$ by
compatible isomorphisms. Then we get the second isomorphism from the first.
\bull

\medskip

\medskip\noindent {\it Proof of Propositions 11 and 12.} Using the isomorphism
$RHK\cong RH\te_{H\cap K} RK$ as $RH\te RK^\circ$-modules (by the evident
maps) whenever $H,K$ are subgroups of a finite group $G$, one has
\smallskip\noindent{$e_IRP_I wRP_{I_0} e_{I_0}\cong RL_I wRL_{I_0}\cong
RL_I w\te_{L_{I^w\cap {I_0}}}RL_{I_0}=$}\hfill\break\rightline{$RP_Ie_I we_{I^w\cap
{I_0}}\te_{P_{I^w\cap {I_0}}}e_{I^w\cap {I_0}}e_{I_0}P_{I_0}.$}

Applying Proposition~10, we get $e_IRP_I\dot
wRP_{I_0}e_{I_0}\cong RP_Ie_{I\cap^w{I_0}} we_{I^w\cap
{I_0}}\te_{P_{I^w\cap
{I_0}}}e_{I^w\cap
{I_0}}RP_{I_0}$ by the map of $RP_I\te RP_{I_0}^\circ$-modules sending
$e_Ix\dot w ye_{I_0}$ to
$xe_{I\cap^w{I_0}}\dot we_{I^w\cap
{I_0}}\te e_{I^w\cap
{I_0}}y$ for $x\in RP_I$, $y\in RP_{I_0}$ and $\dot w\in N$ such that $\dot wT=w$.
Note that our map does not depend on the choice of $\dot w$ inside the class $w\in
N/T$.

On applying the functor $RGe_I\te_{P_I}-$, we get $RGe_I\te_{P_I}e_IRP_I
wRP_{I_0}e_{I_0}\cong RGe_{I\cap^w{I_0}} we_{I^w\cap
{I_0}}\te_{P_{I^w\cap
{I_0}}}e_{I^w\cap
{I_0}}RP_{I_0}$ by the map $xe_I\te e_Iy\dot wze_{I_0}\mapsto
xye_{I\cap^w{I_0}}\dot we_{I^w\cap {I_0}}\te_{P_{I^w\cap
{I_0}}}e_{I^w\cap
{I_0}}z$ for $x\in RG$, $y\in RP_I$, $z\in RP_{I_0}$. On the other hand,
$RGe_{I\cap^w{I_0}}\dot we_{I^w\cap {I_0}}=RG e_{I^w\cap
{I_0}}$ by Dipper-Du-Howlett-Lehrer's theorem of independence (see
[HL]~2.4, [CaEn]~3.10). So the map of Proposition~11 is indeed  defined and an
isomorphism as announced.

Assume now the hypotheses of Proposition~12. Using the restriction map of
$X(L_{I_0})$, it is clear that the map announced in the first statement of
Proposition~12 is the map $RGe_{I_0}\te_{L_{I_0}}\varphi_{J^{w'}\cap {I_0},I^w\cap
{I_0}}$ on
$RGe_{I_0}\te_{L_{I_0}}X(L_{I_0})$, where $\varphi$ denotes the restriction maps of
$X(L_{I_0})$ as a coefficient system on $2^{I_0}$.

To verify the second statement, we assume $W_Jw'=W_Jw$ with $w'\in D_{J{I_0}}$
(unique). It suffices to check that the following square is commutative~:
\smallskip\centerline{$\matrix{X_{I,w}&\to&Y_{I,w}\cr
\downarrow&&\downarrow\cr X_{J,w'}&\to&Y_{J,w'} } $}\smallskip
\noindent  where
horizontal arrows are the isomorphism of Proposition~11, the first vertical arrow
is the restriction map $(X(G)e_{I_0})^I\to (X(G)e_{I_0})^J$ (which actually sends the
term
$X_{I,w}$ in the term $X_{J,w'}$) and the second vertical arrow is the one we have
seen above, that is $RGe_{I_0}\te\varphi_{J^{w'}\cap {I_0},I^w\cap {I_0}}$.

We check the images on a $RG\te RP_{I_0}^\circ$-generator of $X_{I,w}$, namely
 $e_I\te e_I\dot we_{I_0}$.
Going right then down, one gets $e_I\te_{P_I}e_I\dot we_{I_0}\mapsto
e_{I\cap^w{I_0}}\dot we_{I^w\cap {I_0}}\te_{P_{I^w\cap {I_0}}}e_{I^w\cap {I_0}}\mapsto
e_{I\cap ^w{I_0}}\dot we_{I^w\cap {I_0}}\te_{P_{J^{w'}\cap {I_0}}}e_{I^w\cap {I_0}}$.
Going down then right yields $e_I\te_{P_I}e_I\dot we_{I_0}\mapsto
e_I\te_{P_J}e_I\dot we_{I_0} =e_I\te_{P_J}e_I\dot v\dot w'e_{I_0}=e_I\dot
v\te_{P_J}\dot w'e_{I_0}\mapsto e_I\dot ve_{J\cap ^{w'}{I_0}}\dot w'e_{J ^{w'}\cap
{I_0}}\te_{P_{J ^{w'}\cap {I_0}}}e_{I ^{w}\cap {I_0}}$ where $w=vw'$ with $v\in W_J$ such that
$w'$ is $J$-reduced and one defines $\dot v=\dot w(\dot w')^{-1}\in L_J$.

But $e_{J \cap ^{w'} {I_0}}\dot w'e_{J ^{w'}\cap {I_0}}=e_J\dot w'e_{I_0}=e_J\dot w'e_{J
^{w'}\cap {I_0}}$ (Proposition~10 above), so that $e_I\dot ve_{J \cap^{w'} {I_0}}\dot
w'e_{J ^{w'}\cap {I_0}}=e_I\dot ve_J\dot w'e_{J ^{w'}\cap {I_0}}=e_Ie_J\dot we_{J ^{w'}\cap
{I_0}}=e_I\dot we_{J ^{w'}\cap {I_0}}$ so the last term on the second composition (down
then right) is  $e_I\dot w\te_{P_{J ^{w'}\cap {I_0}}}e_{I ^{w}\cap {I_0}}$. This is what we
expect, again by Proposition~10. \bull

\bigskip
\centerline {\bf 3. Alvis-Curtis duality in the
homotopy category.}

\medskip

The main result of this section shows that the derived equivalence of [CaRi]
actually holds in the homotopy category ([O]~Th.~1).

Let us recall two general results that will be useful.
Let $A$ be a ring and $X,Y\in \Cb(A)$, $X_0\in \Cb(A^\circ )$. We refer to
[CaEn]~\S~A1 for the usual notations $X[i]\in \Cb(A)$ ($i\in \Z$), $\Homgr _A(X,Y),
X_0\te_AX\in \Cb(\Z )$. If $A$ is an $R$-algebra for $R$ a
commutative ring, one denotes $X^\vee :=\Homgr_R(X,R[0])$ as an object of
$\Cb(A^\circ )$ in the usual way.

The same proof as [Ri1]~1.1.(a) gives the following

\medskip\noindent {\bf Lemma~14. }  {\sl  Let $C$ be in
$\Cb(A)$ such that 
its homology is concentrated in degree 0 (i.e. $H(C)\cong H^0(C)[0]$ in $\Cb(A)$) and, for
any term
$C'$ of
$C$, both
$\Homgr_{A}(C',C)$ and $\Homgr_{A}(C,C')$ have their homology
concentrated in degree 0.

 Then 
$C\cong H^0(C)[0]$ in $\Kb(A)$.}

\medskip For the strong adjunction below, see [Ro]~\S~2.2.6, [Ri2]~9.2.5. Note that the
isomorphisms can be made explicit (see [Ro]~\S~2.2, [O]~\S~5.1).

\medskip\noindent {\bf Proposition~15. } {\sl Assume $A$, $B$ and $\Lambda$ are
$R$-free $R$-algebras of finite rank. Assume $A$ and $B$ are symmetric (see
[CaEn]~1.19, [Ri2]~9.2.1, [Ro]~2.2.3). All modules are assumed to be $R$-free of finite
rank.

Let $M$ be a bounded complex of $A\te B^\circ $-modules projective on each side.
Let
$N_1$, $N_2$ be objects of $\Cb(B\te \Lambda^\circ )$,  $\Cb(A\te \Lambda^\circ )$
respectively.
Then 

\centerline{${\rm Homgr}_{A\te \Lambda^\circ
}(M\te_BN_1,N_2)\cong {\rm Homgr}_{B\te \Lambda^\circ }(N_1,M^\vee
\te_AN_2).$}}
\medskip

\medskip\noindent {\bf Theorem~16. } (Okuyama, [O]~Th.~1)  {\sl Let $G$ be a finite group
endowed with a split BN-pair of characteristic $p$, let $R =\Z [p^{-1}]$. Recall $X(G)$ in
$\Cb(RG\te RG^\circ )$ (see Notation~8). 
Then 
$$X(G)\te_{RG}X(G)^\vee\cong X(G)^\vee\te_{RG} X(G)\cong RG\ \ {\rm in}\ \ 
\Kb(RG\te RG^\circ )$$ and therefore
$X(G)\te_{RG}-$ induces a splendid equivalence
$\Kb(RG)\to \Kb(RG)$ in the sense of [Ri1].}
\medskip

\medskip\noindent {\it Proof. } We have the claimed
isomorphisms in
$\Db (RG\te RG^\circ)$ (see [CaRi]~5.1 or [CaEn]~4.18). Note that in both references, the
main argument is the following fact~: \hfill\break\noindent {\bf (F)} {\it if $I\inn S$
and
$N$ is a cuspidal
$RL_I$-module, then
$M:=RGe_I\te_{L_I}N\in \mod(RG)$ satisfies $X(G)\te_G M\cong M[-|I|]$ in $\Db
(RG)$ (and even in $\Kb (RG)$).} \hfill\break\noindent
Here, {\it cuspidal} means that, considering $N$ as
a $U_I$-trivial
$RP_I$-module, one has
$e_JN=0$ for any $J\subset I$, $J\not=I$. The rest of the proof of [CaRi]~5.1, essentially
consists in reducing to $R$ being a field, then use the fact that simple
$RG$-modules are quotients of those $M$'s.

 A homotopy
equivalence is preserved by any additive functor, so we may apply
$-\te_{L_I}N\colon \mod(RG\te RL_I^\circ )\to \mod(RG)$ to the first
isomorphism of Theorem~9, and get 
$X(G)\te_GM\cong RGe_I\te_{L_I}X(L_I)\te_{L_I}N$ in $\Kb (RG)$. By
cuspidality of
$N$, one has clearly $X(L_I)\te_{L_I}N=N[-|I|]$. Whence (F).

 Let us now prove
the isomorphisms of Theorem~16 in
$\Kb (RG\te RG^\circ)$. We use induction on $|S|$, the case $G=T$ being trivial.

For the first isomorphism $X(G)\te_{G} X(G)^\vee\cong RG$, in view of
Lemma~14, we have to check that, for any direct summand $C$ of a term of $
X(G)\te_{G} X(G)^\vee $, the $\Homgr_{RG\te RG^\circ }$'s between $C$ and $
X(G)\te_{G} X(G)^\vee$ have their homology in degree 0 only.

The terms of $X(G)\te_{G} X(G)^\vee $ are of the type
$RGe_I\te_{P_I}e_IRG\te_GRGe_J\te_{P_J}e_JRG$ for $I,J\inn S$. This rearranges as
$RGe_I\te_{P_I}e_IRGe_J\te_{P_J}e_JRG =\oplus_{w\in D_{IJ}}
RGe_I\te_{P_I}e_IRP_IwRP_Je_J\te_{P_J}e_JRG$. Applying Proposition~11, one finds
a sum of modules $Y_{I,w}\te_{P_J}e_JRG = RGe_{I^w\cap J}\te_{P_{I^w\cap
J}}e_{I^w\cap J}RG$, each isomorphic to some
$X(G)^{I_0}=E_{I_0}\te_{L_{I_0}}E_{I_0}^\vee$  where $E_{I_0}=RGe_{I_0}$ for
${I_0}\inn S$. 

So we have to check that, for all ${I }\inn S$, both $ \Homgr_{G\times
G^\circ }(E_{I }\te_{L_{I}}E_{I }^\vee ,X(G)\te_{G} X(G)^\vee) $ and
$\Homgr_{G\times G^\circ }( X(G)\te_{G} X(G)^\vee,E_{I }\te_{L_{I }}E_{I }^\vee)
$ have homology in degree 0 only. By Proposition~15, it suffices to check
${\rm Endgr}_{L_I\times G^\circ}(E_{I }^\vee\te_G X(G))$. The bi-additivity of the functor
$\Homgr_A(-,-)$ along with Theorem~9 give ${\rm Endgr}_{L_I\times G^\circ}(E_{I
}^\vee\te_G X(G))\sim {\rm Endgr}_{L_I\times G^\circ}(X(L_I)\te_{L_I}E_{I }^\vee
)$ where $\sim$ denotes homotopy
equivalence in $\Cb (\Z )$

Case $I \not=S$. Our induction hypothesis tells us that $X(L_{I })^\vee\te_{L_{I }}X(L_{I })$
is homotopically equivalent to $RL_{I } [0]$ in
$\Cb (R L_{I }\te RL_{I }^\circ )$. Then Proposition~15 gives ${\rm Endgr}_{L_I\times
G^\circ}(X(L_I)\te_{L_I}E_{I }^\vee )\cong \Homgr_{L_I\times
G^\circ}(X(L_I)^\vee \te_{L_I}X(L_I)\te_{L_I}E_{I }^\vee ,E_{I }^\vee )\sim \Homgr_{L_I\times
G^\circ}(E_{I }^\vee ,E_{I }^\vee )$ which is in degree zero.

Case $I =S$.  One has to check the homology of ${\rm Endgr}_{G\times G^\circ} (X(G))$.
There is a spectral sequence $_IE^{pq}_i$ for the double complex $\Hom_{G\times
G^\circ}(X(G),X(G))$, such that $E_0^{pq}=\Hom_{G\times G^\circ}(X(G)^{-p},X(G)^q)$ and
$E_1^{pq}=H^q(\Homgr_{G\times G^\circ}(X(G)^{-p},X(G)))$ (see for instance [B]~\S~3.4 on
spectral sequences of double complexes). 

Let $J\inn S$. Theorem~9 again gives
$\Homgr_{G\times G^\circ}(X(G)^J,X(G))\cong \Homgr_{G\times
L_J^\circ}(E_J,X(G)\te_G E_J)\sim \Homgr_{G\times
L_J^\circ}(E_J,E_J\te_{L_J}X(L_J))$ . This has $q$-th homology $=0$ for $q>|J|$ since
$X(L_J)^q=0$ for those $q$. Then
$E_1^{pq}=0$ for
$p+q>0$. Whence $H^i({\rm Endgr}_{G\times G^\circ} (X(G)))=0$ for $i>0$. Negative $i$'s are
taken care by the spectral sequence $_{II}E_i^{pq}$ (satisfying 
$E_1^{pq}=H^p(\Homgr_{G\times G^\circ}(X(G),X(G)^q))$).

This completes the proof of $X(G)\te_{G} X(G)^\vee\cong RG$ in $\Kb(RG\te
RG^\circ )$. As for the isomorphism $X(G)^\vee\te_{G} X(G)\cong RG$, it suffices
to do the same with
$X(G)^\vee$ instead of $X(G)$. Note that $X(G)^\vee$ has the same terms as $X(G)$. The
only non trivial fact that we need is a version of the second isomorphism of Theorem~9 for
$X(G)^\vee$. This in turn is a
consequence of the first on applying the functor $M\to
M^\vee$.\bull

\bigskip
\centerline {\bf 4. Hecke algebras.}

\medskip
Let $(W,S)$ be a finite Coxeter group. Let $R$
be a commutative ring. Let $(q_s)_{s\in S}\in (R^\times)^S$ be a family of
invertible elements of $R$ such that $q_s=q_t$ whenever $s,t\in S$ are
$W$-conjugate.

Recall the definition of the Hecke algebra $\H = \oplus_{w\in W}Rh_w$ (see for
instance [GP]~4.4.6) with multiplication obeying the rules

$h_wh_{w'}=h_{ww'}$ when $w,w'\in W$ and lengths add (therefore $h_1=1_\H$),

$(h_s)^2=(q_s-1)h_s+q_s$ when $s\in S$. 

Note that $\H$ is symmetric for the linear
form giving the coordinate on $h_1$ in the above basis (see [GP]~8.1.1).
 When $I\inn S$, $\H_I= \oplus_{w\in
W_I}Rh_w$ is a subalgebra of $\H$ and is also the Hecke algebra associated to
$(W_I,I)$ and same coefficients $q_s$.

Following [LS], one defines the complex $X(\H )$ of $\H$-bimodules, from the
following coefficient system on $2^S$. For $I\inn J\inn S$, let
$X(\H )^I=\H\te_{\H_I}\H$ with restriction maps $\varphi_{JI}\colon X(\H )^I\to
X(\H)^J$ defined by
$\varphi_{JI}(h\te_{\H_I} h') =h\te_{\H_J}h'$.

\medskip\noindent {\bf Theorem~17.} (Okuyama, [O]~4.1)  {\sl   Let ${I_0}\inn S$. The
restriction of $X(\H)$ to
$\H\te_R(\H_{I_0})^\circ $ is isomorphic in $\Kb(\H\te(\H_{I_0})^\circ  )$ with
$\H \te_{\H_{I_0}}X(\H_{I_0})$.}

\medskip\noindent {Proof. } The proof is very similar to the one of
Theorem~9. Let us abbreviate $\te_{\H_I}$ as $\te_I$ in what follows.
Since
$h_w$'s multiply as elements of
$W$ when lengths add, one has $\H =\oplus_{w\in D_{I{I_0}}}\H_Ih_w\H_{I_0}$ as
$\H_I,\H_{I_0}$-bimodule.
Then the restriction to $\H\te (\H_{I_0})^\circ $ of $X(\H )^I$ is
$\H\te_{I}\oplus_{w\in D_{I{I_0}}}\H_Ih_w\H_{I_0}\cong \oplus_{w\in D_{I{I_0}}}\H
.(1\te_{ I}h_w).\H_{I_0}$.

Corresponding to Propositions~11 and 12, one has the following.

Let $I\inn S$, $w\in D_{I{I_0}}$.
Denote $X(\H )_{I,w}=\H\te_{ {I}}{\H_{I}}h_w{\H_{{I_0}}}$.

Then $X(\H )_{I,w}\cong \H\te_{ {I^w\cap {I_0}}}\H_{I_0}$ by a map sending
$1\te_{ {I}}h_w$ to $h_w\te_{ {I^w\cap {I_0}}}1$. To see that, note first that if
$v\in W_{I^w\cap {I_0}}$, then ${}^wv\in W_{I}$ and therefore
$h_wh_v=h_{wv}=h_{{}^wv}h_w$. The claimed isomorphism then corresponds to
the composition of the following (explicit) isomorphisms
$\H
\te_{ I}\H_Ih_w\H_{I_0} = \H
\te_{ I}\H_I\te_{ {I\cap^w{I_0}}}\H_{I\cap^w{I_0}}h_w\H_{I_0}  =\H
\te_{ I}\H_Ih_w\te_{ {I^w\cap {I_0}}}\H_{I_0} =\H
h_w\te_{ {I^w\cap {I_0}}}\H_{I_0}=\H
\te_{ {I^w\cap {I_0}}}\H_{I_0}$.

If $I^w\cap {I_0}\inn J^{w'}\cap {I_0}$ with $J\inn S$ and $w'\in D_{J{I_0}}$, we have the
restriction map
$\H\te_{ {I^w\cap {I_0}}}\H_{I_0}\to
\H\te_{ {J^{w'}\cap {I_0}}}\H_{I_0}$.

If moreover $W_Jw =W_Jw'$, the above corresponds to the restriction
map of $X(\H )$, $\varphi_{JI}\colon X(\H )_{I,w}\to X(\H )_{J,w'}$ through the above
isomorphism.
To check this, one evaluates at $1\te h_w$. Denote $v\in W_J$ such that $w=vw'$
with lengths adding. One composition (isomorphism, then restriction map) gives
$1\te_{ {I}} h_w\mapsto h_w\te_{ {I^w\cap {I_0}}} 1\mapsto
h_w\te_{ {J^{w'}\cap {I_0}}} 1$, the other composition gives
$1\te_{ {I}} h_w\mapsto 1\te_{ {J}} h_w=1\te_{ {J}}
h_vh_{w'}=h_v\te_{ {J}} h_{w'}\mapsto h_vh_{w'}\te_{ {J^{w'}\cap {I_0}}} 1$. 

One may now replace $_\H X(\H )_{ {\H_{I_0}}}$ by the complex associated to the
coefficient system on $2^S$ associating to $I$ the bi-module
$\oplus_w\H\te_{I^w\cap I_0}\H_{I_0}$, a sum over $D_{II_0}$, with restriction maps
defined by all the maps
$\H\te_{ {I^w\cap {I_0}}}\H_{I_0}\to
\H\te_{ {J^{w'}\cap {I_0}}}\H_{I_0}$ for $I\inn J$, $w'\in D_{JI_0}$ and $W_Jw =W_Jw'$
(which implies ${I^w\cap {I_0}}\inn J^{w'}\cap {I_0}$).

Then $_\H X(\H )_{ {\H_{I_0}}}\to\!\!\!\!\to \H\te_{ {{I_0}}} X(\H_{I_0})$ by a map
sending
$\H\te_{ {I^w\cap {I_0}}}\H_{I_0}$ to $0$ if $W_Iw\not\inn W_{I_0}$, and
$\H\te_{ {I}}\H_{I_0}\to \H\te_{ {{I_0}}}\H_{I_0}\te_{ {I}}\H_{I_0}$ the evident (onto)
map when $I\inn I_0$. The image is the complex associated to the coefficient
system
$\H\te_{I_0} X(\H_{I_0})$ on
$2^{I_0}$.

Using now the notation of \S~1, the kernel corresponds to the complex associated
with the coefficient system on $2^S$, $I\mapsto
Z^I=\oplus_{b\in
\AP ,S(b)=I}\H\te_{I_0(b)}\H_{I_0} $ with restriction maps sending
$\H\te_{I_0(b)}\H_{I_0}$ into $\H\te_{I_0(b')}\H_{I_0}$ by
$\H\te_{I_0}\varphi_{I_0(b'),I_0(b)}$ only when $S(b')\nni S(b)$ and $b'=b\cup S(b')$. 
The contractibility of
$Z$ is immediate by Theorem~4. Then Proposition~13 implies our claim.
\bull

\medskip\noindent {\bf Remark~18. }  The above theorem allows to give a proof of
the main result in [LS] and [PS]~\S~3.1 in the
following way.

Taking ${I_0}=\vid$ in Theorem~17, one gets that $X(\H )$ has its homology concentrated
in degree 0. So this homology is the
kernel of
$\ddd ^0\colon \H\te_R\H\to \oplus_{s\in S}\H\te_{\H_s}\H$, that is the
intersection of the kernels of the restriction maps $\varphi_{ s , \vid }\colon
\H\te_R\H\to
\H\te_{\H_s}\H$ defined by $h\te h'\mapsto h\te_{\H_s}h'$. 

An element $x\in \H\te_R\H$ writes in a unique way $x=\sum_{w\in W}h_w\te
x_w$ with $x_w$'s in $\H$. 
For $s\in S$, using the partition $W= D_{\vid
,s}\sqcup D_{\vid ,s}s$ and the law of $\H$, it is clear that $\varphi_{  s,\vid } (x)=0$
if and only if
$x_{ws}=-(h_s)^{-1}x_w$ for any $w\in D_{\vid , s}$. Using a reduced decomposition
of each $w\in W$, one then gets that
$x\in
\ker(\ddd ^0)$ if and only if $x_w=(-1)^{l(w)}(h_w)^{-1}x_1$ for any $w\in W$.

Define
$\xi :=
\sum_{w\in W}(-1)^{l(w)}h_w\te (h_w)^{-1}\in
\H\te_R\H$. We now have
$\ker(\ddd ^0)=\xi.\H\cong
\H$ as right $\H$-module. 

Using the partition
$W=D_{s,\vid}\sqcup sD_{s,\vid}$ and the formula for $(h_s)^2$, it is clear that
$h_s \xi h_{s}=-q_s \xi $, for any $s\in S$. So we get
\hfill\break\centerline{$H^0(X(\H))=\ker(\ddd ^0)=\H .\xi =\xi .\H\cong
{}_\alpha\H$} where
$\alpha$ is the automorphism of
$\H$ sending $h_s$ to $-q_s(h_s)^{-1}$. 

The $R$-dual $X(\H )^\vee$ has similar properties, so $H(X(\H ))\te_\H
H(X(\H )^\vee )\cong {}_\H\H_\H [0]$. On the other hand, the bi-projectivity of the terms of
$X(\H)$ implies that the natural map $H(X(\H )\te_\H X(\H )^\vee )\to H(X(\H ))\te_\H
H(X(\H )^\vee )={}_\H\H_\H [0]$ is an isomorphism (see [B]~3.4.4).

{\sl }
\medskip\noindent {\bf Theorem~19. } (Okuyama, [O]~Th.~2)  {\sl 
$X(\H)\te_\H -$ induces an
auto-equivalence of the homotopy category $\Kb(\H)$.}

\medskip\noindent {\it Proof. } The proof is similar to the one of Theorem~16
with $\H$ replacing $RG$, ${}_\H\H_{\H_{I_0}}$ replacing $RGe_{I_0}$,
Theorem~17 replacing Theorem~9, and Remark~18
giving the equivalence in the derived category.
Also
$\H$ is symmetric, which allows to use Proposition~15.

Note that we need to have identical statements for $X(\H)$ and $X(\H)^\vee$. As in
the end of proof of Theorem~16, this is done by applying the contravariant
functor $M\mapsto M^\vee$ and using the isomorphism $X(\H)\cong
X(\H_{I_0})\te_{\H_{I_0}}\H$ in $\Kb(\H_{I_0}\te\H^\circ )$. This in turn is deduced from
Theorem~17 by noting
that there is a (covariant) functor
$M\mapsto M^\iota$ from $\mod(\H)$ to $\mod(\H^\circ )$ corresponding to the
isomorphism $\H\cong \H^\circ $ (as $R$-algebras) which is defined by $h_s\mapsto
h_s$ (and therefore $h_w\mapsto h_{w^{-1}}$).\bull
\bigskip

\centerline {\bf 5. Generalized Steinberg module and the Solomon-Tits theorem.}
\medskip

The context (and notations) are now again the ones of \S 2 and 3, where $G$ is a finite
group with a split BN-pair. Define {\rm St}$(G) $ as the
object of
$\Cb (\Z G)$ associated to the coefficient system on
$2^S$ defined by ${\rm
St}(G)^I=\Z G/P_I$ and $\varphi_{JI}(gP_I)=gP_J$ for $I\inn J\inn S$ and
$g\in G$.

Compare the following with an old lemma on the Steinberg character ([DM]~9.2).

\medskip\noindent {\bf Theorem~20.} (Okuyama, [O]~3.7)  {\sl Let $I\inn S$. Then ${\rm
Res}^G_{P_I}{\rm St} (G) \cong {\rm Ind}^{P_I}_{L_I}{\rm
St} (L_I)$ in $\Kb (\Z P_I)$.
}

\medskip\noindent {\it Remark.} Note that, for $I=\vid$ (and therefore
$P_I=B$, $L_I=T$), the theorem implies that St$(G)$ has homology only in
degree 0 (a theorem of Solomon-Tits, see [CuRe]~66.33). Note also that the proof below
simplifies a lot when $I=\vid$. 

\medskip\noindent {\it Proof of Theorem~20.} We prove the statement with $I_0$ instead
of $I$.
 We actually check a right module version of (i) with 
${\rm
St} (G)^I=\Z (P_I\bs G)$, noting that the left module version follows by the same type
of considerations as in the end of proofs of Theorem~9 and  Theorem~16.

Note also that the definition of St$(G)$ can be made using any system of
parabolic subgroups containing a given Borel subgroup $B^{g_0}$. Denote
$B^-=B^{w_S}$,
$P_I^-=B^-W_IB^- =U_I^-.L_I$ where $U^-_I=U^{w_S}\cap U^{w_Sw_I}$. 

One may
assume that
${\rm
St} (G)$ (resp. ${\rm
St} (L_{I_0})$) is associated to the coefficient system on $2^S$
(resp. 
$2^{I_0}$) defined by $I\mapsto [P_I^-]\Z  G$, (resp. $I\mapsto
[ L_{I_0}\cap P_I^-] \Z  L_{I_0}$) where $[F]$ denotes the sum of elements
of $F$ in
$\Z G$ for any subset $F\inn G$.

Let us define a surjective map ${\rm Res}_{P_{I_0}}^G {\rm
St} (G)\to
{\rm
St} (L_{I_0})\te_{L_{I_0}}\Z P_{I_0}$ in
$\Cb (\Z P_{I_0}{}^\circ )$, split in each degree and with contractible kernel. Then our
claim will follow by Proposition~13.

Noting that $P^-_I\bs G/P_{I_0}=W_I\bs W/W_{I_0}$ is in bijection with $D_{II_0}$, one
has ${\rm Res}_{P_{I_0}}^G({\rm
St} (G)^I)= \oplus_{w\in D_{II_0}}[P_I^-]w\Z P_{I_0}$ with
connecting map ${\rm Res}_{P_{I_0}}^G({\rm
St} (G)^I\to {\rm
St} (G)^J)$ sending $[P_I^-]wg$ to
$[P_J^-]wg$ for any $g\in {P_{I_0}} $ and $I\inn J\inn S$.

Let us define a map $\pi\colon {\rm Res}_{P_{I_0}}^G {\rm
St} (G)\to
{\rm
St} (L_{I_0})\te_{L_{I_0}}\Z P_{I_0}$ by sending $[P_I^-]w\Z P_{I_0}$ to $0$ except
when
$I\inn I_0$ and $w=1$ in which case we use the isomorphism $\pi_I\colon
[P_I^-]\Z P_{I_0}\cong [P_I^-\cap L_{I_0}]\Z P_{I_0} = [P_I^-\cap L_{I_0}]\Z L_{I_0}\te_{
L_{I_0}}\Z P_{I_0}$ due to the fact that $P_I^-\cap P_{I_0} = P_I^-\cap L_{I_0}$ (apply for
instance [CaEn]~Exercise~2.4, or see proof of the Lemma below). This is clearly a
surjective
$\Z P_{I_0}^\circ$-homomorphism on coefficients, split at each $I$. In order to
show that it commutes with connecting maps, we have to check for any $I\inn J\inn
S$, the equality $\pi_J\varphi_{JI}=\varphi_{JI}\pi_I$ where the first $\varphi$ is the
connecting map in St$(G)$ and the second is the one in
St$(L_{I_0})\te_{L_{I_0}}\Z P_{I_0}$ (seen as a coefficient system on $2^S$ by
extending trivially on $2^S\sm 2^{I_0}$). If
$J\not\inn I_0$, both sides are 0. If $I\inn J\inn I_0$, one gets $[P_I^-]g\mapsto
[P_J^-]g\mapsto [P_J^-\cap L_{I_0}]g $ and $[P_I^-]g\mapsto [P_I^-\cap L_{I_0}]g\mapsto
[P_J^-\cap L_{I_0}]g $ if $g\in P_{I_0}$,
$[P_I^-]g\mapsto [P_J^-]g\mapsto 0$ and $[P_I^-]g\mapsto 0\mapsto 0$ if $g\in
G\sm P_{I_0}$.    

The kernel $Z$ of our map is a graded sum of $\Z P_{I_0}^\circ$-modules
$Z_{b}=[P_I^-]w\Z P_{I_0}$ for $b=W_Iw\in \AP$, $Z_b$ being at degree $|I|$. To show that
$Z$ is contractible, one imitates the proof of Theorem~4. Recall
$\ss\in {\rm End}_\Z(\Z\AP )$ from Proposition~3, which is expressed by $\ss
(b)=\sum_{b'\in\AP}m_{bb'}b'$ with
$m_{bb'}\in\Z$ for $b,b'\in\AP$. 

Fix $b=W_Iw$ and $b'=W_{I'}w'$ in $\AP$ with $I,I'\inn S$ and $w\in D_{II_0}$, $w'\in D_{I'I_0}$
such that $m_{bb'}\not=0$ . Proposition~3 and Remark~7 tell us that
$I_0\cap I^w\inn I_0\cap (I')^{w'}$ and
$ w_Sw w_{I_0}
\leq_{\bf r} w_Sw' w_{I_0}$. We use that through the following consequence

\medskip\noindent {\bf Lemma.} {\sl $(P_I^-)^{w}\cap P_{I_0}\inn (P_{I'}^-)^{w'}\cap P_{I_0}
$.}
\medskip
This allows to define $\varphi_{b'b}\colon Z_b = [P_{I}^-]w\Z P_{I_0}\to Z_{b'}=
[P_{I'}^-]w'\Z P_{I_0}$ as the only $(\Z P_{I_0})^\circ$-homomorphism sending $[P_{I}^-]w$
to
$[P_{I'}^-]w'$. This behaves like restriction maps~:
$\varphi_{b''b'}\circ \varphi_{b'b}=\varphi_{b''b}$ whenever $b,b',b''\in\AP$ with
$m_{bb'}\not=0$ and
$m_{b'b''}\not=0$. Defining now
$\bar\ss\colon Z\to Z$ as
$\sum_{b'\in\AP}m_{bb'}\varphi_{b'b}$ on $Z_b$, the same proof as for Theorem~4 shows
that the equations ${\bf
E5}_b$ satisfied by the $m_{bb'}$'s imply
$\bar\ss\partial +\partial\bar\ss = \Id$ on $Z$.\bull

\medskip\noindent {\it Proof of the Lemma.} Replacing $I$ and $I'$ by their $w_S$-conjugate
and $w$, $w'$ by $w_Sww_{I_0}$, $w_Sw'w_{I_0}$ respectively, we have to check that
$(P_I)^{w}\cap P_{I_0}\inn (P_{I'})^{w'}\cap P_{I_0} $ as soon as $w\in D_{II_0}$, $w'\in
D_{I'I_0}$ satisfy $I^w\cap I_0\inn (I')^{w'}\cap I_0$ and $w'\leq_{\bf r} w$. 

Applying twice [CaEn]~2.27.(i), one gets $L_{I^w\cap I_0}\inn P_I^w\cap P_{I_0}\inn
P_{I^w\cap I_0} =L_{I^w\cap I_0}.U$ and
$^wU\cap P_I\inn U$. Therefore $P_I^w\cap
P_{I_0}=L_{I^w\cap I_0}.( U^w\cap U)$. Now $I^w\cap I_0\inn (I')^{w'}\cap I_0$,
while $w'\leq_{\bf r} w$ implies $  U^w\cap U\inn
 U^{w'}\cap U$ (apply [CaEn]~2.3.(iii), 2.23.(i)).
\bull \bigskip

\bigskip
\centerline {\bf References.}
\smallskip
[B] D. Benson, {\it Representations and Cohomology II: Cohomology of Groups and
Modules}, Cambridge Univ. Press, Cambridge, 1991.

[CaEn] M. Cabanes and  M. Enguehard, {\it Representation Theory of Finite Reductive
Groups}, Cambridge Univ. Press, Cambridge, 2004.

[CaRi] {M. Cabanes and J. Rickard}, Alvis--Curtis duality as
an equivalence of derived categories, in {\sl Modular Representation Theory
of Finite Groups}, (M.J. Collins, B.J. Parshall, L.L. Scott, eds), de Gruyter,
(2001), 157--174.


[CuRe] C.W. Curtis and  I. Reiner, {\it Methods of Representation Theory with Applications to
Finite Groups and Orders}, Wiley, 1987.

[DM] F. Digne and J. Michel, {\it Representations of finite groups of Lie type}, Cambridge,
1991.

{[GP]} M. Geck and G. Pfeiffer, {\it Characters of Finite Coxeter Groups and
Iwahori-Hecke Algebras}, Oxford, 2000.

[HL] {B. Howlett and G. Lehrer}, On Harish-Chandra
induction for modules of Levi subgroups, {\sl J. Algebra}, {\bf 165}, (1994),
172--183.

[LS] M. Linckelmann and S. Schroll, A two-sided $q$-analogue of the Coxeter
complex, {\sl J. Algebra}, {\bf 289}, 1 (2005), 128-134.

[O] T. Okuyama, On conjectures on complexes of some module
categories related to Coxeter complexes, preprint, August 2006, 25 p.

[PS] B. Parshall and L. Scott, Quantum Weyl reciprocity for cohomology, {\sl Proc.
London Math. Soc.} (3), {\bf 90}, (2005), 655-688.

[Ri1] J. Rickard, Splendid equivalences : derived categories and permutation
modules, {\sl Proc. London Math. Soc.} (3), {\bf 72}, (1996), 331-358.

[Ri2] J. Rickard, Triangulated categories in the modular representation theory of
finite groups, {\it in} S. K\"onig and A. Zimmermann, {\it Derived Equivalences for
Group Rings}, {LNM}, {\bf 1685}, Springer,  1998, 177-198

[Ro] {R. Rouquier}, Block theory via stable and Rickard
equivalences, in {\it Modular Representation
Theory of Finite Groups}, (M.J. Collins, B.J. Parshall, L.L.
Scott, eds), Walter de Gruyter, (2001),  101--146.

\footline={\hfil}

\bigskip{}\bigskip{}\bigskip{}\bigskip{}\bigskip{}\bigskip{}\bigskip{}\bigskip{}

\noindent{Marc CABANES}, {UNIVERSITE DENIS DIDEROT}

\noindent{175, rue du Chevaleret, F-75013 PARIS, FRANCE}

\noindent\it E-mail adress : {\tt cabanes@math.jussieu.fr}
\end